\documentclass[10pt]{article} 
\usepackage{amsmath, amssymb, amsthm, mathrsfs}
\usepackage{enumerate}
\usepackage{colordvi}
\usepackage{tikz}
\usepackage{graphicx}
\usepackage{hyperref}
\newtheorem{claim}{Claim}[section]
\newtheorem{theorem}[claim]{Theorem}
\newtheorem{lemma}[claim]{Lemma}

\newtheorem{corollary}[claim]{Corollary} 
\usepackage{color}

\definecolor{Myred}{cmyk}{0.0,1.0,1.0,0.00}
\definecolor{Mypurple}{rgb}{0.5,0.0,0.5}

\newtheorem{definition}{Definition}

\begin{document}
%\begin{center}
%{\Large{\textbf{Spectral convergence of Neumann Laplacian in the perforated domain}}}
\title{Magnetic Neumann Laplacian on a domain with hole}
%\bigskip
\author
{
Diana Barseghyan$^{1,2}$\footnote{corresponding author}, Swanhild Bernstein$^{3}$ and Baruch Schneider$^{1}$
}
\date{\small $^{1}$ Department of Mathematics, University of Ostrava,  30.dubna 22, Ostrava 70103, Czech Republic\\ \small $^{2}$ Department of Theoretical Physics\\
Nuclear Physics Institute, Czech Academy of Sciences,
25068 \v{R}e\v{z}, Czech Republic\\ \small $^{3}$ Institute of Applied Analysis, TU Bergakademie Freiberg, Akademiestrasse 6, Freiberg 09599, Germany\\
E-mails:\, diana.barseghyan@osu.cz, Swanhild.Bernstein@math.tu-freiberg.de,  baruch.schneider@osu.cz
}
\maketitle

\begin{abstract}
This article gives a domain with a small compact set of removed and the magnetic Neumann Laplacian on such set. The main theorem of this article shows the description of the holes which do not change the spectrum drastically. In this article we prove that the spectrum of the magnetic Neumann Laplacian converges in the Hausdorff distance sense to the spectrum of the original operator defined on the unperturbed domain.

\end{abstract}
\bigskip


%%%%%%%%%%%%%%%%%%%%%%

\section{Introduction}
\label{s:intro}

Let $\Omega\subset \mathbb{R}^2$ be a bounded domain and a compact set $K\subset \Omega$. We denote 
$\Omega_K:= \Omega\setminus K$.  
The magnetic Neumann Laplacian $$H_{\Omega_K}(A_{\mathrm{magn}})=(i \nabla+A_{\mathrm{magn}})^2$$ is defined via the quadratic form
\[
\int_{\Omega_K} |i\nabla u+A_{\mathrm{magn}}u|^2dxdy,\quad u\in \mathcal{H}^1(\Omega_K),
\]
where $A_{\mathrm{magn}}$ is a magnetic vector potential generating the magnetic field 
$$\mathrm{rot}(A_{\mathrm{magn}})(x, y)= B(x, y).$$

However for applications, a better do not use different notation for the magnetic potentials corresponding to the  restrictions of the same magnetic field. The first motivation arises for studying the magnetic Laplacian from the mathematical theory of superconductivity, see e.g \cite{Saint}. The second motivation is to understand at which point there is an analogy between the a magnetic Laplacian and a non-magnetic Laplacian. In this paper we will focus on problems with magnetic Laplacian.

It is a common expectation that small perturbations of the physical situation lead only to a small change of the spectrum. In the case of domain perturbations this is largely true for Dirichlet boundary conditions while the Neumann case is more delicate.

In Neumann case even small perturbations may cause abrupt change of the spectrum. For example, such an effect is observed when the hole has a "split-ring" geometry \cite{S15}. The split ring (even being very small) may produce additional eigenvalues having nothing in common with the eigenvalues of the Neumann Laplacian on the unperturbed domain.

The problems with the Neumann obstacles having more general geometry appeared in \cite{CP20}. The authors required the hole to satisfy the so called ”uniform extension property”. This is a requirement, which means that $\mathcal{H}^1$- functions on the domain with a hole can be extended to $\mathcal{H}^1$ function on the unperturbed domain and the norm of this extension operator does not depend on the hole diameter. In this case the authors established the spectral convergence.

In \cite{BSH22} it was considered the spectral problem for the Neumann Laplacian on the domain with hole with zero Lebesgue measure. Obviously such domain does not satisfy the uniform extension property but under some additional assumption on hole one is able to prove the spectral convergence.

In the present work we discuss the similar situation where one is not able to guarantee the validity of the uniform extension property.  Here we do not restrict ourselves that the hole must have zero Lebesgue measure. 
Moreover, we extend the class of operators: instead of the Laplacian we consider the magnetic Laplace operators.

Before the description of our results let as also mention paper \cite{BM18}, which makes use imposing Neumann boundary conditions on the boundary of hole but in higher dimension.  The situation with a smooth hole is considered in \cite{BE22}.

The second section of the present work is devoted to the description of the geometry of holes which do not change the spectrum drastically. More precisely, we will be interested in the spectral properties of the magnetic Neumann Laplacian defined on a two-dimensional bounded domain with a single hole $K_\varepsilon$ (for a fixed  parameter $\varepsilon$) having "standard" geometry. In this section we present the results about the spectral convergence of the magnetic Neumann Laplacian on $\Omega_{K_\varepsilon}$ as $\varepsilon\to 0$ in terms of the Hausdorff distance under some additional assumptions on the geometry of $K_\varepsilon$.
In the third section we present the main tools of the spectral convergence of operators on varying Hilbert spaces. In the fourth and fifth sections we prove our results.

\section{Main results}
\setcounter{equation}{0}
\bigskip

Since we are working with the geometry of domains, we introduce the additional \textbf{property$^*$}: for any $(x_0, y_0)\in \Omega_K$
at least one of the following conditions takes place

\begin{itemize}
\item Let $l(x_0)=\{x= x_0\}$. Then at least one of the half-lines $l(x_0)\cap\{y\ge y_0\}$ and $l(x_0)\cap\{y\le y_0\}$ has no intersection with $K$.   

\item Let $h(y_0)=\{y= y_0\}$. Then at least one of the half-lines $h(y_0)\cap\{x\ge x_0\}$ and $h(y_0)\cap\{x\le x_0\}$ has no intersection with $K$.
\end{itemize}
{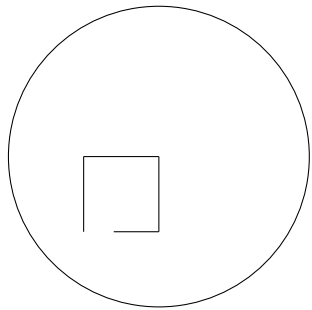
\begin{figure}
\centering
\begin{tikzpicture}
\draw (2,2) -- (2,1);
\draw (1,2) -- (2,2);
\draw (1,2) -- (1,1);
\draw (1.4,1) -- (2,1);
\draw (2,2) circle (2cm);
\end{tikzpicture}
\caption{Example of $\Omega_K$ not satisfying \textbf{property$^*$}}\label{figure1}
\end{figure}}

%%%%%%%%%%%%\textbf{PICTURE}
\bigskip

We collect in the following theorem some important results.

\begin{theorem}\label{main}
Let $\Omega$ be an open bounded domain in $\mathbb{R}^2$ and let $p\in \Omega$ be some fixed point. Suppose that ${\mathbb B}_\varepsilon\subset \Omega$ is a ball with center at $p$ and radius $\varepsilon>0$. Let $K= K(\varepsilon)\subset \mathbb B_\varepsilon$ be a compact set. Moreover, suppose that $\Omega_K$ satisfies {\textbf{property$^*$}}. Let $H_\Omega(A_{\mathrm{magn}})$ and $H_{\Omega_{K_\varepsilon}}(A_{\mathrm{magn}})$ be the magnetic Neumann Laplacians defined on $\Omega$ and $\Omega_{K_\varepsilon}$, respectively. We assume that on some subset of $\Omega$ containing ${\mathbb B}_\varepsilon$ the vector potential $A_{\mathrm{magn}}$ is twice continiously differentiable.
Then for small enough $\varepsilon$ there exists $\eta(\varepsilon)>0$ with $\eta(\varepsilon)\to 0$ as $\varepsilon\to 0$ such that the following spectral convergence takes place
$$\overline{d} \left(\sigma_{\bullet} \left(H_{\Omega_{K_\varepsilon}}(A_{\mathrm{magn}})\right), \, \sigma_{\bullet} \left(H_\Omega(A_{\mathrm{magn}})\right)\right) \le \eta(\varepsilon),$$
where $\overline{d}$ is defined in (\ref{Hausdorff2}) and $\sigma_\bullet (\cdot)$ denotes either the entire spectrum, the essential or the discrete spectrum. Moreover, the multiplicity of the discrete spectrum is preserved.
\end{theorem}
\qed
\\

\noindent The previous result motivates the following consequences:

%As consequences of Theorem \ref{main}:

\begin{corollary}
Suppose that $H_\Omega(A_{\mathrm{magn}})$ has purely discrete spectrum denoted by $\lambda_k(\Omega)$ (repeated according to multiplicity). Then the infimum of the essential spectrum of $H_{\Omega_{K_\varepsilon}}(A_{\mathrm{magn}})$ tends to infinity and there exists $\eta_k(\varepsilon)>0$ with $\eta_k(\varepsilon)\to 0$ as $\varepsilon\to 0$ such that
$$|\lambda_k(\Omega)- \lambda_k(\Omega_{K_\varepsilon})|\le \eta_k(\varepsilon)$$
for small enough $\varepsilon$. Here, $\lambda_k(\Omega_{K_\varepsilon})$ denotes the discrete spectrum of $H_{\Omega_{K_\varepsilon}}(A_{\mathrm{magn}})$ (below the essential spectrum) repeated according to multiplicity.
\end{corollary}

\begin{corollary} The Hausdorff distance between the spectra of $H_{\Omega_{K_\varepsilon}}(A_{\mathrm{magn}})$ and $H_\Omega(A_{\mathrm{magn}})$ converges to zero on any compact interval $[0, \Lambda]$. 
\end{corollary}

The proof of Theorem \ref{main} is based on Theorem \ref{POSTtm} and on the following theorem:

\begin{theorem}\label{Second}
Under the assumptions of Theorem \ref{main}
for small enough $\varepsilon$ the operators $H_\Omega(A_{\mathrm{magn}})$ and $H_{\Omega_{K_\varepsilon}}(A_{\mathrm{magn}})$ are $ \mathcal{O}\left(\varepsilon^{1/6}\right)$  close of order $2$.
\end{theorem}

\section{Main tool of the spectral convergence  of  operators on varying Hilbert spaces}
\setcounter{equation}{0}

We begin this section by reviewing some basic facts which assures a spectral convergence for two operators having the different domains. For more information we refer the reader to  \cite{P06}. 

To a Hilbert space $H$ with inner product $(\cdot, \cdot)$ and norm $\|\cdot\|$ together with a non-negative, unbounded, operator $A$, we associate the scale of Hilbert spaces
\[
H_k:=\mathrm{Dom}((A+1)^{k/2}),\quad \|u\|_k:=\|(A+ 1)^{k/2}u\|,\,\,k\ge 0,
\]
where $1$ is the identity operator.

We think of $(H', A')$ being some perturbation of $(H, A)$ and want to lessen the assumption such that the spectral properties are not the same but still are close.

\begin{definition}(see \cite{P06})\label{Post} 
Suppose we have linear operators
\begin{gather*} J: H \longrightarrow H', \quad\quad J_1: H_1\longrightarrow H_1'\\
J': H'\longrightarrow H,\quad\quad J_1': H_1' \longrightarrow H_1. \end{gather*}
Let $\delta> 0$ and $k\ge 1$.
We say that $(H, A)$ and $(H', A')$ are $\delta$-close of order $k$ iff the following conditions are fulfilled:

\begin{gather}\label{1} \|J f- J_1 f\|_0\le \delta \|f\|_1, \\\label{2}  |(J f, u)-(f, J' u)|\le \delta \|f\|_0 \|u\|_0,\\\label{3}  \|u-J J' u\|_0\le \delta \|u\|_1,\\\label{4} \quad \|Jf\|_0\leq 2\|f\|_0,\,\,\, \|J'u\|_0\le 2\|u\|_0,
\\\label{5'} \|(f- J'J f)\|_0\le \delta \|f\|_1,\\
\label{6} \| J' u- J_1'u\|_0\le \delta \|u\|_1\\\label{7}  |a(f, J_1' u)- a'(J_1 f, u)|\le \delta \|f\|_k \|u\|_1,\end{gather}
for all $f, u$ in the appropriate spaces. Here, $a$ and $a'$ denote the sesquilinear forms associated to $A$ and $A'$.
\end{definition}

We denote by $d_{\mathrm{Haussdorff}}(A, B)$ the Hausdorff distance for subsets $A, B\subset\mathbb{R}$
\begin{equation}\label{Haus.}
d_{\mathrm{Haussdorff}}(A, B):=\mathrm{max} \left\{\underset{a\in A} {\mathrm{sup}}\,d(a, B),\,\underset{b \in B} {\mathrm{sup}}\, d(b, A)\right\},
\end{equation}
where $d(a, B):= \mathrm{inf}_{b\in B} |a - b|$. We set
\begin{equation}\label{Hausdorff2}\overline{d} (A, B):= d_{\mathrm{Hausdorff}}\left((A+1)^{-1}, (B+1)^{-1}\right)\end{equation}
for closed subsets of $[0, \infty)$. For the next result, which originates with the  work of O. Post \cite{P06} we have the following spectral convergence theorem in terms of the distance $\overline{d}$.
\begin{theorem}\cite{P06}\label{POSTtm}
There exists $\eta(\delta)>0$ with $\eta(\delta)\to 0$ as $\delta\to 0$ such that
\begin{equation}\label{Post}\overline{d}(\sigma_\bullet(A), \, \sigma_\bullet (A'))\le \eta(\delta)\end{equation}
for all pairs of non-negative operators and Hilbert spaces $(H, A)$ and $(H', A')$ which are $\delta$-close. Here, $\sigma_\bullet (A)$ denotes either the entire spectrum, the essential or the discrete spectrum of $A$. Furthermore, the multiplicity of the discrete spectrum, 
$\sigma_{\mathrm{disc}}$, is preserved, i.e. if $\lambda\in \sigma_{\mathrm{disc}}(A)$ has multiplicity $\mu>0$, then there exist $\mu$ eigenvalues (not necessarily all distinct) of operator $A'$ belonging to interval $(\lambda-\eta(\delta), \lambda+\eta(\delta))$.
\end{theorem}
\qed

We now turn to the proof of Theorem \ref{Second}.

\section{Proof of Theorem\,\ref{Second}}

We split the proof of this statement into number of steps.
\bigskip
\\
Step \,(i).\,{\it Construction of the mappings $J, J', J_1, J_1'$}.

\bigskip

It is easy to notice that
$H= L^2(\Omega),\,H'= L^2(\Omega_{K_\varepsilon}),\,A= A'=(i \nabla+A_{\mathrm{magn}})^2$;
$H_1,\,H_1'$ correspond to Sobolev spaces $\mathcal{H}^1(\Omega)$ and $\mathcal{H}^1(\Omega_{K_\varepsilon})$ and
$H_2= \mathrm{Dom}(H_\Omega(A_{\mathrm{magn}}))$. The norm $\|\cdot\|_0$ corresponds with the
$L^2$ norm and
$$\|u\|_1=(\|u\|^2_0+\|i \nabla u+A_{\mathrm{magn}}u\|_0^2)^{1/2},\quad
\|f\|_2=(\|(i \nabla+A_{\mathrm{magn}})^2f+f\|_0^2)^{1/2}.$$ 
%%%%%%%%%

We define $J u=J_1 u=u|_{\Omega_{K_\varepsilon}}$ for all $u\in H$ and  
{$J' u =u \chi_{\Omega_{K_\varepsilon}}$ for all $u\in H'$.}

Let us now construct the mapping $J_1': H_1' \to H_1$.
Without loss of generality, assume that the ball $\mathbb B_\varepsilon$ mentioned in Theorem \ref{main} and Theorem \ref{Second} is centered at the origin. Let $\epsilon\in (\varepsilon, 2\varepsilon)$ be a number to be chosen later and let $\mathbb B_\epsilon \supset \mathbb B_\varepsilon$ be the ball with center again at the origin and radius 
$\epsilon$, $\Omega_\epsilon:=\Omega\backslash \mathbb B_\epsilon$.

We are going to construct mapping $J_1'$ first for smooth functions. For any $v\in  C^\infty(\Omega_{K_\varepsilon})$  we define

\begin{equation*}J_1' v:=\begin{cases} v,\quad \text{on} \quad \Omega_\epsilon, \\ \frac{r}{\epsilon} \tilde{v}(\epsilon, \varphi), \quad \text{on} \quad \mathbb B_\epsilon,\end{cases}\end{equation*}
where $\tilde{v}(r, \varphi)=v(r \cos \varphi, r \sin \varphi)$.

Now let us construct the mapping $J_1' u$ for any $u\in H_1'$.  Employing the approximation method described in \cite[Thm.2, 5.3.2]{E10}, for the fixed sequence $\{\eta_k\}_{k=1}^\infty$ converging to zero we construct the sequence 
$v_{\eta_k}\in C^\infty(\Omega_{K_\varepsilon})$ which satisfies
\begin{equation}\label{converging}
\|u-v_{\eta_k}\|_1< \eta_k\,\|u\|_1.\end{equation}

Let us mention that in view of the inequalities (\ref{J1'est.}) and (\ref{eps.}) which will be proved later it follows  for any smooth function $v$
\begin{eqnarray*}\|J_1'v\|_1^2= \int_\Omega |i\nabla J_1'v+A_{\mathrm{magn}}J_1'v|^2\,d x\,d y+\int_\Omega |J_1'v|^2\,d x\,d y\\
\le 2\int_\Omega |\nabla J_1'v|^2\,d xd y+2\||A_{\mathrm{magn}}|\|_{L^\infty(\Omega)}^2\int_\Omega |J_1'v|^2\,d x\,d y\le \overline{C}(\varepsilon) \|v\|_1^2,\end{eqnarray*}
where $\overline{C}(\varepsilon)$ is some constant.

Therefore using the density of space $H_1$ we are able to define
\begin{equation}\label{approximation}
J_1' u= \underset{k\to \infty}{\mathrm{lim}} \,J_1'v_{\eta_k}.
\end{equation}
\\
Step \,(ii).\, \it {The conditions (\ref{1})-(\ref{7}) hold for the mappings $J, J', J_1, J_1'$}.

\bigskip
\rm Indeed, we have that the estimates (\ref{1})-(\ref{4}) are satisfied with $\delta=0$. 
\\\\
We now prove (\ref{5'}), i.e. {\it under the assumptions stated in Theorem \ref{Second} inequality (\ref{5'}) is satisfied with $\delta=\mathcal{O}(\sqrt{\varepsilon})$ for small enough $\varepsilon$.}
\\

In view of our construction we have 
\begin{eqnarray}\nonumber
\|f- J' J f\|_0^2=\int_\Omega |f- J' f\chi_{\Omega_{K_\varepsilon}}|^2\,d x\,d y=\int_\Omega |f- f\chi_{\Omega_{K_\varepsilon}}|^2\,d x\,d y\\\label{start 6}=\int_{K_\varepsilon}|f|^2\,d x\,d y\le \int_{B_\varepsilon}|f|^2\,d x\,d y.
\end{eqnarray}

To proceed further we need to employ the following lemma (the proof will be given in Appendix) applied with 
$\Gamma_\varepsilon=\emptyset$:
\begin{lemma}\label{auxiliary!!}
Suppose the assumptions of Theorem \ref{main}. To avoid the misunderstanding we denote the compact hole inside of $B_\varepsilon$ by $\Gamma_\varepsilon$. Then for any function $u\in\mathcal{H}^1(\mathbb{B}_\epsilon\setminus \Gamma_\varepsilon)$ the following inequality takes place
$$\int_{\mathbb{B}_\epsilon\setminus \Gamma_\varepsilon}|u(x, y)|^2\,d x\,d y\le C_1\varepsilon\int_{\Omega\setminus \Gamma_\varepsilon}(|i\nabla u+A_{\mathrm{magn}}u|^2+|u|^2)\,d x\,d y,$$
where $C_1>0$ is a constant depending on the distance between the boundary of 
$\mathbb{B}_\epsilon$ and the boundary of $\Omega$ and the magnetic potential $A_{\mathrm{magn}}$. 
\end{lemma}

By Lemma \ref{auxiliary!!}  and (\ref{start 6}) we prove (\ref{5'}).
\\

Using (\ref{approximation}) and the completeness of $C^\infty(\Omega_{K_\varepsilon})$ in Sobolev space 
$\mathcal{H}^1(\Omega_{K_\varepsilon})$ it is enough to prove (\ref{6}) for $u\in C^\infty(\Omega_{K_\varepsilon})$.

Taking into account that $J' u=0$ on $K_\varepsilon$ and $J_1' u= u$ on $\Omega\setminus \mathbb{B}_\epsilon$ one has
\begin{gather}
\nonumber
\|J' u- J_1' u\|_0^2= \int_{\Omega_{K_\varepsilon}}|u-J_1'u|^2\,d x\,d y+ \int_{K_\varepsilon}|J_1' u|^2\,d x\,d y\\\nonumber =\int_{\mathbb{B}_\epsilon\setminus K_\varepsilon}|u-J_1'u|^2\,d x\,d y+ \int_{K_\varepsilon}|J_1' u|^2\,d x\,d y
\\\label{ineq}\le 2\int_{\mathbb{B}_\epsilon\setminus K_\varepsilon}|u|^2\,d x\,d y+2\int_{\mathbb{B}_\epsilon}|J_1' u|^2\,d x\,d y+ \int_{K_\epsilon}|J_1' u|^2\,d x\,d y.
\end{gather}

To estimate the first integral of the right-hand side of (\ref{ineq}) we apply Lemma \ref{auxiliary!!} with 
$\Gamma_\varepsilon=K_\varepsilon$.

Let us now investigate the second term in the right-hand side of (\ref{ineq}).
Passing to polar coordinates one gets
\begin{eqnarray*}
\nonumber &&\int_{\mathbb B_\epsilon}|J_1'u|^2\,d x\,d y = \int_0^\epsilon \int_0^{2\pi} r \left|\frac{r}{\epsilon}
\tilde{u}(\epsilon, \varphi)\right|^2\,d r\,d \varphi\\&&\le\int_0^\epsilon \int_0^{2\pi} r |\tilde{u}(\epsilon, \varphi)|^2\,d r\,d \varphi\le \epsilon \int_0^\epsilon \int_0^{2\pi} 
|\tilde{u}(\epsilon, \varphi)|^2\,d r\,d \varphi\\&&= \epsilon^2 \int_0^{2\pi} 
|\tilde{u}(\epsilon, \varphi)|^2\,d \varphi.
\end{eqnarray*}
Taking into account that $\epsilon \int_0^{2\pi} 
|\tilde{u}(\epsilon, \varphi)|^2\, d \varphi$ coincides with the curvilinear integral 
$\int_{\partial \mathbb{B}_\epsilon} |u|_{\partial \mathbb{B}_\epsilon}^2\, d \mu$, where $|u|_{\partial \mathbb{B}_\epsilon}$ is the trace of $|u|$ on the circle $\partial \mathbb{B}_\epsilon$, the above bound performs to
\begin{equation}\label{ugamma}
\int_{\mathbb B_\epsilon}|J_1'u|^2\,d x\,d y\le\epsilon
\int_{\partial \mathbb{B}_\epsilon}|u|_{\partial \mathbb{B}_\epsilon}^2\,d \mu.
\end{equation}

The suitable bound for the above estimate is guaranteed by the following trace inequality \cite{E10}:

\begin{lemma}\label{trace}
Let $v\in \mathcal{H}^1(\Omega\backslash \mathbb{B})$, where $\mathbb{B}$ is a ball with center in zero and radius $\tau>0$.
Then there exists a constant $C_2>0$ depending on the diameter of $\Omega$ and the distance between the boundary of $\mathbb{B}$ and the boundary of $\Omega$ such that
$$\int_{\partial \mathbb B} 
|v|_{\partial \mathbb B}^2\,d \mu\le C_2 \int_{\Omega\backslash \mathbb{B}}(|\nabla v|^2+|v|^2)\,d x\,d y.$$
\end{lemma}

Let us show that for every $v\in\mathcal{H}^1(\omega),\,\omega\subset\mathbb{R}^2$, and a magnetic vector potential $A$ one has
\begin{equation}\label{Comp.}\int_\omega(|\nabla v|^2+|v|^2)\,d x\,d y\le C_3(A)\int_\omega (|i \nabla v+Av|^2+|v|^2)\,d x\,d y,\end{equation}
where $C_3(A)$ is given by (\ref{C3'}).

Indeed, since
$$\int_\omega |i\nabla v+Av|^2\,d x\,d y\ge \frac{1}{2}\int_\omega|\nabla v|^2\,d x\,d y- 2\||A|\|^2_\infty\int_\omega|v|^2\,d x\,d y$$
then 
$$\int_\omega(|\nabla v|^2+|v|^2)\,d x\,d y\le 2\int_\omega|i\nabla v+A v|^2\,d x\,d y+ (4\||A|\|_\infty^2+1)\int_\omega|v|^2\,d x\,d y,$$
which establishes our aim with \begin{equation}\label{C3'}C_3(A)=\mathrm{max}\{2, 4\||A|\|_\infty^2+1\}.\end{equation}

Inequality (\ref{ugamma}) together with Lemma \ref{trace}, (\ref{Comp.}) and  fact that $\epsilon\le 2\varepsilon$ gives 
\begin{equation}\label{J1'est.}
\int_{\mathbb B_\epsilon}|J_1'u|^2\,d x\,d y\le 2C_2C_3\varepsilon\|u\|_1^2,
\end{equation}
where 
\begin{equation}\label{C3}
C_3= C_3(A_{\mathrm{magn}}).
\end{equation}

Since
$$\int_{K_\varepsilon}|J_1'u|^2\,d x\,d y\le\int_{B_\epsilon}|J_1'u|^2\,d x\,d y$$
then 
$$\int_{K_\varepsilon}|J_1'u|^2\,d x\,d y\le 2C_2C_3\varepsilon \|u\|_1^2.$$

Combining the inequality above together with (\ref{J1'est.}) and Lemma \ref{auxiliary!!} the right-hand side of (\ref{ineq}) can be estimated as follows
$$
\|J' u- J_1'u\|_0^2\le 2C_2C_3\varepsilon \|u\|_1^2 
$$
which proves (\ref{6}) with $\delta=\mathcal{O}(\sqrt{\varepsilon})$.\\\\

We now give the proof of the estimate (\ref{7}), i.e. {\it under the assumptions stated in Theorem \ref{Second} inequality (\ref{7}) takes place with $k=2$ and $\delta=\mathcal{O}(\varepsilon^{1/6})$ for small enough 
$\varepsilon$.}
\\
We have
\begin{eqnarray}\nonumber
|a(f, J_1'u)- a'(J_1f, u)|\\\nonumber=\left|\int_\Omega(i\nabla f+A_{{\mathrm{magn}}}f)\overline{(i\nabla (J_1'u)+A_{{\mathrm{magn}}}J_1'u)}\,d x\,d y-\int_{\Omega_{K_\varepsilon}}(i\nabla (J_1f)+A_{{\mathrm{magn}}}J_1f)\overline{(i\nabla u+A_{{\mathrm{magn}}}u)}\,d x\,d y
\right|\\\nonumber=\biggl|\int_{\Omega_{\mathbb B_\epsilon}}(i\nabla f+A_{{\mathrm{magn}}}f)\overline{(i\nabla u+A_{{\mathrm{magn}}}u)}\,d x\,d y+\int_{\mathbb B_\epsilon}(i\nabla f+A_{{\mathrm{magn}}}f)\overline{(i\nabla (J_1'u)+A_{{\mathrm{magn}}}J_1'u)}\,d x\,d y\\\nonumber-\int_{\Omega_{\mathbb B_\epsilon}}(i\nabla f+A_{{\mathrm{magn}}}f)\overline{(i\nabla u+A_{{\mathrm{magn}}}u)}\,d x\,d y-\int_{\mathbb B_\epsilon\backslash K_\varepsilon}(i\nabla f+A_{{\mathrm{magn}}}f)\overline{(i\nabla u+A_{{\mathrm{magn}}}u)}\,d x\,d y\biggr|\\\nonumber=\biggl|\int_{\mathbb B_\epsilon}(i\nabla f+A_{{\mathrm{magn}}}f)\overline{(i\nabla (J_1'u)+A_{{\mathrm{magn}}}J_1'u)}\,d x\,d y-\int_{\mathbb B_\epsilon\backslash K_\varepsilon}
(i\nabla f+A_{{\mathrm{magn}}}f)\overline{(i\nabla u+A_{{\mathrm{magn}}}u)}\,d x\,d y\biggr|\\\label{7'}\le\left|\int_{\mathbb B_\epsilon}(i\nabla f+A_{{\mathrm{magn}}}f)\overline{(i\nabla (J_1'u)+A_{{\mathrm{magn}}}J_1'u)}\,d x\,d y\right|+\left|\int_{\mathbb B_\epsilon\backslash K_\varepsilon}
(i\nabla f+A_{{\mathrm{magn}}}f)\overline{(i\nabla u+A_{{\mathrm{magn}}}u)}\,d x\,d y\right|.\end{eqnarray}

Let us estimate each term of (\ref{7'}). Starting from the first term one gets
\begin{eqnarray}\nonumber\left|\int_{\mathbb B_\epsilon}(i\nabla f+A_{{\mathrm{magn}}}f)\overline{(i\nabla (J_1'u)+A_{{\mathrm{magn}}}J_1'u)}\,d x\,d y\right|\\\label{8}\le\left(\int_{\mathbb B_\epsilon}|i\nabla f+A_{{\mathrm{magn}}}f|^2\,d x\,d y\right)^{1/2}\left(\int_{\mathbb B_\epsilon}|i\nabla (J_1'u)+A_{{\mathrm{magn}}}J_1'u|^2\,d x\,d y\right)^{1/2}.\end{eqnarray}

Dealing with the first term of the right-hand side of the above inequality we have 
\begin{eqnarray}\nonumber\int_{\mathbb B_\epsilon}|i\nabla f+A_{{\mathrm{magn}}}f|^2\,d x\,d y\le2\int_{\mathbb B_\epsilon}|\nabla f|^2\,d x\,d y+2\int_{\mathbb B_\epsilon}|A_{{\mathrm{magn}}}|^2|f|^2\,d x\,d y\\\label{9}\le2\int_{B_\epsilon}|\nabla f|^2\,d x\,d y+2\||A_{{\mathrm{magn}}|}\|^2_\infty\int_{B_\epsilon}|f|^2\,d x\,d y.\end{eqnarray}
 
The first and second terms in (\ref{9}) can be estimated via the following lemmas (the proofs are given in Appendix):
\begin{lemma}\label{magnetic}
For any function $g\in \mathrm{Dom}(H_\Omega(A_{\mathrm{magn}}))$ the following estimate takes place
\[
\int_{\mathbb B_\epsilon}|\nabla g|^2\,d x\,d y\le C_5\epsilon^{4/3}\int_\Omega|(i\nabla+A_{\mathrm{magn}})^2 g+g|^2\,d x\,d y,
\]
with the constant $C_5$ depends on $\Omega$ and magnetic potential $A_{\mathrm{magn}}$.
 \end{lemma}

\begin{lemma}\label{magnetic1}
For any $g\in \mathrm{Dom}(H_\Omega(A_{\mathrm{magn}}))$ the following inequality takes place
\[
\int_{\mathbb B_\epsilon}|g|^2\,d x\,d y\le C_6\epsilon^{4/3}\int_{\Omega}|(i\nabla+A_{\mathrm{magn}})^2 g+g|^2\,d x\,d y,
\]
with the constant $C_6$ depends on $\Omega$ and magnetic potential $A_{\mathrm{magn}}$.
\end{lemma}

This together with the fact that $\epsilon\le 2\varepsilon$ estimates the right-hand side of (\ref{9}) as follows
\begin{equation}\label{rhs9}\int_{\mathbb B_\epsilon}|i\nabla f+A_{{\mathrm{magn}}}f|^2\,d x\,d y\le 2^{4/3}C_7\varepsilon^{4/3}\int_\Omega|(i\nabla+A_{{\mathrm{magn}}})^2 f+f|^2\,d x\,d y=2^{4/3}C_7\varepsilon^{4/3}\|f\|_2^2\end{equation}
with $C_7= 2C_5+2\||A_{{\mathrm{magn}}}|\|_\infty^2C_6$.

To proceed further with the proof of an upper bound of (\ref{8}) we need to estimate integral $\int_{\mathbb B_\epsilon}|i\nabla (J_1'u)+A_{{\mathrm{magn}}}J_1'u|^2\,d x\,d y$. One has
\begin{eqnarray}\nonumber \int_{\mathbb B_\epsilon}|i\nabla (J_1'u)+A_{{\mathrm{magn}}}J_1'u|^2\,d x\,d y\le
2\int_{\mathbb{B}_\epsilon}|\nabla (J_1'u)|^2\,d x\,d y+2\int_{\mathbb{B}_\epsilon}|A_{{\mathrm{magn}}}|^2|J_1'u|^2\,d x\,d y\\\label{10}\le2\int_{\mathbb{B}_\epsilon}|\nabla (J_1'u)|^2\,d x\,d y+2\||A_{{\mathrm{magn}}}|\|_\infty^2\int_{\mathbb{B}_\epsilon}|J_1'u|^2\,d x\,d y.\end{eqnarray}

Passing to polar coordinates in the first integral in the right-hand side of (\ref{10}) one gets
\begin{eqnarray}\nonumber&&\int_{\mathbb B_\epsilon} |\nabla J_1' u|^2\,d x\,d y=
\int_{\mathbb B_\epsilon} \left(\left|\frac{\partial (J_1' u)}{\partial x}\right|^2+ \left|\frac{\partial (J_1' u)}{\partial y}\right|^2 \right)\,d x\,d y\\\nonumber&&
= \int_0^\epsilon \int_0^{2\pi}
r \left|\frac{1}{\epsilon} \tilde{u}(\epsilon, \varphi)\,\cos \varphi -\frac{1}{\epsilon} \frac{\partial \tilde{u}}{\partial \varphi}(\epsilon, \varphi)\, \sin \varphi\right|^2\,d r\,d \varphi\\\nonumber&&+ \int_0^\varepsilon \int_0^{2\pi}
r \left|\frac{1}{\epsilon} \tilde{u}(\epsilon, \varphi)\,\sin \varphi +\frac{1}{\epsilon} \frac{\partial \tilde{u}}{\partial \varphi}(\epsilon, \varphi)\, \cos \varphi\right|^2\,d r\,d \varphi\\\nonumber&& \le \frac{4}{\epsilon} \int_0^\epsilon \int_0^{2\pi} \left(|\tilde{u}(\epsilon, \varphi)|^2+\left|\frac{\partial \tilde{u}}{\partial \varphi}(\epsilon, \varphi)\right|^2\right)\,d r\,d \varphi\\\label{J 1'}&& =4 \int_0^{2\pi}|\tilde{u}(\epsilon, \varphi)|^2d\varphi+ 4 \int_0^{2\pi}\left|\frac{\partial \tilde{u}}{\partial \varphi}(\epsilon, \varphi)\right|^2\,d \varphi.
\end{eqnarray}

As in proof of (\ref{6}) we note that
$$4 \int_0^{2\pi}|\tilde{u}(\epsilon, \varphi)|^2\,d \varphi=\frac{4}{\epsilon}\,\int_{\partial \mathbb B_\epsilon} 
|u|_{\partial \mathbb B_\epsilon}^2\,d \mu,$$
which can be estimated using Lemma \ref{trace} and (\ref{Comp.}) as follows:

\begin{equation}\label{first term}
4 \int_0^{2\pi}|\tilde{u}(\epsilon, \varphi)|^2\,d \varphi\le \frac{4}{\epsilon}C_4C_3\|u\|_1^2.
\end{equation}

\bigskip
To obtain the upper bound for the second term of \ref{J 1'} we use the following auxiliary lemma (the proof is given in {Appendix}):
\begin{lemma}\label{auxiliary1}
Let $\mathbb B_{2\varepsilon}$ and $\mathbb B_\varepsilon$ be the balls with center in origin and radii 
$\varepsilon$ and $2\varepsilon$. Let $g\in \mathcal{H}^1(\mathbb B_{2\varepsilon}\setminus \mathbb B_\varepsilon)$. Then there exists $\tau\in (\varepsilon, 2\varepsilon)$ such that
$$\int_0^{2\pi} \left|\frac{\partial \tilde{g}}{\partial \varphi}(\tau, \varphi)\right|^2\,d \varphi \le 4 \int_{\mathbb B_{2\varepsilon}\setminus \mathbb B_\varepsilon}(|\nabla g|^2+|g|^2)\,d x\,d y,$$
where $\tilde{g}(r, \varphi):= g(r\cos\varphi, r\sin\varphi)$.
\end{lemma}

Now we are able to choose the number $\epsilon$. Let us apply Lemma \ref{auxiliary1} for $g=u$. In case if mentioned in above lemma $\tau$ belongs to interval 
$(\varepsilon, 3\varepsilon/2]$ then we choose $\epsilon$ as the supremum of all such numbers in
$(\varepsilon, 3\varepsilon/2]$. In the opposite case if $\tau\in (3\varepsilon/2, 2\varepsilon)$ then let $\epsilon$ be the infimum of such numbers. Since $u$ is smooth function then inequality (\ref{auxiliary1}) is satisfied with $\tau=\epsilon$. 

Combining this together with (\ref{Comp.}) one gets
\begin{equation}\label{sec.term.}\int_0^{2\pi}\left|\frac{\partial \tilde{u}}{\partial \varphi}(\epsilon, \varphi)\right|^2\,d \varphi\le4C_3 \|u\|_1^2.\end{equation}
 
Hence, by virtue of (\ref{J 1'}), (\ref{first term}) and (\ref{sec.term.})
\begin{equation}\label{eps.}\int_{\mathbb B_\epsilon} |\nabla J_1' u|^2\,d x\,d y\le\frac{4}{\epsilon}(C_4C_3+ 16C_3)\|u\|_1^2.\end{equation}

In view of (\ref{J1'est.}), (\ref{10}), (\ref{eps.}) and the fact that $\epsilon\ge\varepsilon/2$ we arrive
\begin{equation}\label{11}
\int_{\mathbb B_\epsilon}|i\nabla (J_1'u)+A_{{\mathrm{magn}}}J_1'u|^2\,d x\,d y\le
\left(\frac{8}{\epsilon}(C_4C_3+ 16C_3)
+4\||A_{{\mathrm{magn}}}|\|_\infty^2C_4C_3\varepsilon\right)\|u\|_1^2.\end{equation}

Finally, employing (\ref{rhs9}) and the above estimate in inequality (\ref{8}) it follows that
\begin{equation}
\nonumber\left|\int_{\mathbb B_\epsilon}(i\nabla f+A_{{\mathrm{magn}}}f)\overline{(i\nabla (J_1'u)+A_{{\mathrm{magn}}}J_1'u)}\,d x\,d y\right|\\\label{f.term}\le2^{13/6}(C_7C_3(C_4+16))^{1/2}\varepsilon^{1/6} \|f\|_2\|u\|_1. 
\end{equation}

Let us now consider the second term in (\ref{7'}). By virtue of (\ref{rhs9}) one has
\begin{eqnarray*}\nonumber\left|\int_{\mathbb B_\epsilon\backslash K_\varepsilon}
(i\nabla f+A_{{\mathrm{magn}}}f)\overline{(i\nabla u+A_{{\mathrm{magn}}}u)}\,d x\,d y\right|\\\nonumber\le\left(\int_{\mathbb B_\epsilon}
|i\nabla f+A_{{\mathrm{magn}}}f|^2dxdx\right)^{1/2}\left(\int_{\mathbb B_\epsilon\backslash K_\varepsilon}|i\nabla u+A_{{\mathrm{magn}}}u|^2\,d x\,d y\right)^{1/2}\\\le 2^{2/3}(C_5)^{1/2}\varepsilon^{2/3}\|f\|_2\|u\|_1.
\end{eqnarray*}

Hence the right-hand side of inequality (\ref{7}) for small enough $\varepsilon$ satisfies
$$
\mathrm{r.h.s.}(\ref{7})= \mathcal{O}\left(\varepsilon^{1/6}\right)\,\|f\|_2\, \|u\|_1,
$$
which ends the proof.\qed
%%%%%%%%

\section{Appendix}

In this section we give the proofs of Lemmas \ref{auxiliary!!}, \ref{magnetic}-\ref{auxiliary1}. Since the outline of our proofs are a slight modification of the proofs in \cite[Appendix]{BSH22}, therefore we refer \cite{BSH22} for identical parts.
\subsection{Proof of Lemma \ref{auxiliary!!}}
\setcounter{equation}{0}

We are going to use the following result:

\begin{lemma}\label{omega.}[The proof will be given in the end of this section]
Let $\Omega_{\Gamma_\varepsilon}$ be an open domain satisfying \textbf{property$^*$}.
Then for any $u\in \mathcal{H}^1(\Omega_{\Gamma_\varepsilon})$ and for almost any
 $(x_0, y_0)\in \Omega_{\Gamma_\varepsilon}$ the following statement takes place: 
\begin{eqnarray}\label{curve1}
|u(x_0, y_0)|\le C_8 \left(\int_{l(x_0)\setminus \Gamma_\varepsilon}\left|\frac{\partial u}{\partial y}(x_0, z)\right|^2\,d z+ \int_{l(x_0)\setminus \Gamma_\varepsilon} |u(x_0, z)|^2\,d z\right)^{1/2},\\\nonumber  \text{if} \quad (x_0, y_0) \quad\text{satisfies first condition of \textbf{property$^*$}} 
\\\label{curve2}
|u(x_0, y_0)|\le C_8 \left(\int_{h(y_0)\setminus \Gamma_\varepsilon}\left|\frac{\partial u}{\partial z}(z, y_0)\right|^2+ \int_{h(y_0)\setminus \Gamma_\varepsilon} |u(z, y_0)|^2\,d z\right)^{1/2}, \\\nonumber\text{if} \quad (x_0, y_0) \quad\text{satisfies second condition of \textbf{property$^*$}},
\end{eqnarray}
where $C_8>0$ is a constant depending on the diameter of $\Omega$ and the distance between point $(x_0, y_0)$ and the boundary of 
$\Omega$.
\end{lemma}

In view of the mentioned result for any $(x, y)\in \mathbb{B}_\epsilon\backslash \Gamma_\varepsilon$
\begin{eqnarray*}|u(x, y)|^2\le C_8^2 \left(\int_{l(x)\setminus \Gamma_\varepsilon}\left|\frac{\partial u}{\partial y}(x, z)\right|^2\,d z+ \int_{l(x)\setminus \Gamma_\varepsilon} |u(x, z)|^2\,d z\right)\\+C_8^2 \left( \int_{h(y)\setminus \Gamma_\varepsilon}\left|\frac{\partial u}{\partial z}(z, y)\right|^2\,d z+ \int_{h(y)\setminus \Gamma_\varepsilon} |u(z, y)|^2\,d z\right).
\end{eqnarray*}

Therefore we obtain
\begin{eqnarray*}
\int_{\mathbb{B}_\epsilon\backslash \Gamma_\varepsilon}|u(x, y)|^2\,d x\,d y\le C_8^2 \int_{\mathbb{B}_\epsilon\backslash \Gamma_\varepsilon}\int_{l(x)\setminus \Gamma_\varepsilon}\left(\left|\frac{\partial u}{\partial z}(x, z)\right|^2+|u(x, z)|^2\right)\,d z\,d x\,d y\\+C_8^2 \int_{\mathbb{B}_\epsilon\backslash \Gamma_\varepsilon} \int_{h(y)\setminus \Gamma_\varepsilon}\left(\left|\frac{\partial u}{\partial z}(z, y)\right|^2+|u(z, y)|^2\right)\,d z\,d x\,d y\\
\le C_8^2\int_{-\epsilon}^\epsilon\int_{\{y: (x, y)\in \mathbb{B}_\epsilon\backslash \Gamma_\varepsilon\}}\int_{l(x)\setminus \Gamma_\varepsilon}\left(\left|\frac{\partial u}{\partial z}(x, z)\right|^2+|u(x, z)|^2\right)\,d z\,d x\,d y\\+ C_8^2 \int_{-\epsilon}^\epsilon\int_{\{x: (x, y)\in \mathbb{B}_\epsilon\backslash \Gamma_\varepsilon\}}\int_{h(y)\setminus \Gamma_\varepsilon}\left(\left|\frac{\partial u}{\partial z}(z, y)\right|^2+|u(z, y)|^2\right)\,d z\,d y\,d x\\\le
2C_8^2\epsilon  \int_{\Omega\backslash \Gamma_\varepsilon}(|\nabla u|^2+ 2|u|^2)\,d x\,d y\le 4C_8^2\epsilon \int_{\Omega\backslash \Gamma_\varepsilon}(|\nabla u|^2+ |u|^2)\,d x\,d y.\end{eqnarray*} 

This together with (\ref{Comp.}) and the fact that $\epsilon\le2\varepsilon$ proves the lemma.

\qed

\subsection{Proof of Lemma \ref{magnetic}}
\setcounter{equation}{0}

To proceed with a proof we need the following auxiliary material  \cite[Lemma\,4.9]{MK06}:
\begin{lemma}\label{marchenko}
Let $\Pi\subset \mathbb{R}^n$ be a convex set and let $G$ and $Q$ be arbitrary measurable sets in $\Pi$ with $\mu\, (G)\ne 0$.
Then, for all $v\in \mathcal{H}^1(\Pi)$, the following inequality holds:
\begin{eqnarray}\label{Marchenko'}
&&\int_Q |v|^2\,d x\,d y\\ \nonumber
&&\le \frac{2\mu\, (Q)}{\mu \,(G)} \int_G |v|^2\,d x\,d y
+\frac{C(n) (d(\Pi))^{n+1} (\mu\, (Q))^{1/n}}{\mu \,(G)} \int_\Pi |\nabla v|^2\, d x\,d y,
\end{eqnarray}  
where $d(\Pi)$ is the parameter of $\Pi$, $\mu$ is the Lebesque measure on $\mathbb{R}^n$, and the constant $C(n)$ depends only on the dimension of $\mathbb{R}^n$.
\end{lemma}

Let $\Pi\subset \mathbb{R}^n$ be a convex set and let $G$ and $Q$ be arbitrary measurable sets in $\Pi$ with $\mu\, (G)\ne 0$.
Applying Lemma(\ref{marchenko}) for $Q=\mathbb B_\epsilon$ and $G=\Pi= G_\epsilon$ we have
\begin{eqnarray*}&& \int_{\mathbb B_\epsilon} \left|\frac{\partial g}{\partial x}\right|^2 \,d x\,d y\le 
\frac{2 \mu \,(\mathbb B_\epsilon)}{\mu\, (G_\epsilon)} \int_{G_\epsilon} \left|\frac{\partial g}{\partial x}\right|^2\,d x\,d y\\&&
+\frac{C(2) (d\,(G_\epsilon))^3 (\mu\, (\mathbb B_\epsilon))^{1/2}}{\mu \,(G_\epsilon)} \int_{G_\epsilon}\left|\nabla\left(\frac{\partial g}{\partial x}\right)\right|^2\,d x\,d y
\\
&&=\frac{2 \mu \,(\mathbb B_\epsilon)}{\mu\, (G_\epsilon)} \int_{G_\epsilon} \left|\frac{\partial g}{\partial x}\right|^2\,d x\,d y+
 \frac{C(2) (d\,(G_\epsilon))^3 (\mu\, (\mathbb B_\epsilon))^{1/2}}{\mu \,(G_\epsilon)} \int_{G_\epsilon}\left(\left|\frac{\partial^2 g}{\partial x^2}\right|^2+ \left|\frac{\partial^2 g}{\partial x \partial y}\right|^2\right)\,d x\,d y\end{eqnarray*}
 and
 \begin{eqnarray*}&&
 \int_{\mathbb B_\epsilon} \left|\frac{\partial g}{\partial y}\right|^2 \,d x\,d y \le \frac{2 \mu \,(\mathbb B_\epsilon)}{\mu\, (G_\epsilon)} \int_{G_\epsilon} \left|\frac{\partial g}{\partial y}\right|^2\,d x\,d y\\&&+
 \frac{C(2) (d\,(G_\epsilon))^3 (\mu\, (\mathbb B_\epsilon))^{1/2}}{\mu \,(G_\epsilon)} \int_{G_\epsilon}\left|\nabla\left(\frac{\partial g}{\partial y}\right)\right|^2\,d x\,d y
\\&&=\frac{2 \mu \,(\mathbb B_\epsilon)}{\mu\, (G_\epsilon)} \int_{G_\epsilon} \left|\frac{\partial g}{\partial y}\right|^2\,d x\,d y+
 \frac{C(2) (d\,(G_\epsilon))^3 (\mu\, (\mathbb B_\epsilon))^{1/2}}{\mu \,(G_\epsilon)} \int_{G_\epsilon}\left(\left|\frac{\partial^2 g}{\partial y \partial x}\right|^2+ \left|\frac{\partial^2 g}{\partial y^2}\right|^2\right)\,d x\,d y.
\end{eqnarray*}

Combining the above inequalities we arrive

\begin{eqnarray*}\nonumber&&
\int_{\mathbb B_\epsilon}|\nabla g|^2 \,d x\,d y\le \frac{2 \mu \,(\mathbb B_\epsilon)}{\mu\, (G_\epsilon)} \int_{G_\epsilon} |\nabla g|^2\,d x\,d y \\\nonumber&&+
 \frac{C(2) (d\,(G_\epsilon))^3 (\mu\, (\mathbb B_\epsilon))^{1/2}}{\mu \,(G_\epsilon)} \int_{G_\epsilon}\left(\left|\frac{\partial^2 g}{\partial x^2}\right|^2+ 2\left|\frac{\partial^2 g}{\partial x \partial y}\right|^2+ \left|\frac{\partial^2 f}{\partial y^2}\right|^2\right)\,d x\,d y\\\nonumber&& \le  \frac{(\mu\,(\mathbb B_\epsilon))^{1/2}}{\mu\, (G_\epsilon)}\left(2(\mu\,(\mathbb B_\epsilon))^{1/2}+ C(2) (d\,(G_\epsilon))^3 \right)\|g\|_{\mathcal{H}^2(G_\epsilon)}^2
 \\
&&= \frac{\sqrt{\pi} \epsilon }{\mu\, (G_\epsilon)}\left(2 \sqrt{\pi} \epsilon+ C(2) (d\,(G_\epsilon))^3 \right)\|g\|_{\mathcal{H}^2(G_\epsilon)}^2.
\end{eqnarray*} 

Let us choose $G_\epsilon=\{r:\, 0\le r\le\epsilon^{1-\alpha}\}$ with some $\alpha\in \left(\frac{1}{2}, 1\right)$ to be chosen later. One can easily notice that $\mathbb B_\epsilon \subset G_\epsilon$. The the above inequality becomes

\begin{eqnarray}\nonumber
\int_{\mathbb B_\epsilon}|\nabla g|^2\,d x\,d y\le 2\epsilon 
\left(\epsilon^{2\alpha-1}+ \frac{4 C(2)}{\sqrt{\pi}} \epsilon^{1-\alpha}\right)\, \|g\|^2_{\mathcal{H}^2(G_\epsilon)}\\\label{Gepsilon}\le 2\mathrm{max}\left\{1, \frac{4 C(2)}{\sqrt{\pi}}\right\}\epsilon\left(\epsilon^{2\alpha-1}+ \epsilon^{1-\alpha}\right)\, \|g\|^2_{\mathcal{H}^2(G_\epsilon)}.
\end{eqnarray}

Let us investigate the function $F(\alpha):= \epsilon^{2\alpha}+ \epsilon^{2-\alpha}$ on interval $\left(\frac{1}{2}, 1\right)$. It attains its minimum at $\alpha_0= \frac{2}{3}- \frac{1}{3}\frac{\ln 2}{\ln \epsilon}$ and takes the value $F(\alpha_0)=\epsilon^{4/3}\left(\frac{1}{4^{1/3}}+ 2^{1/3}\right)$. 

Let us now choose $\alpha= \alpha_0$ in inequality (\ref{Gepsilon}). Then 
\begin{equation}\label{epsilon1}\int_{\mathbb B_\epsilon}|\nabla g|^2\,d x\,d y\le 2\epsilon^{4/3} \left(\frac{1}{4^{1/3}}+ 2^{1/3}\right)\mathrm{max}\left\{1, \frac{4 C(2)}{\sqrt{\pi}}\right\}\,\|g\|_{\mathcal{H}^2(G_\epsilon)}^2.\end{equation}
 
We employ the interior regularity theorem \cite{B16} :

\begin{theorem}(Interior Regularity Theorem.)
Suppose that $h\in \mathcal{H}^1(\Omega)$ is a weak solution to $L h=w$ where $L$ is the elliptic operator on $\Omega$ with $C^1$ smooth coefficients. Suppose that $w\in L^2(\Omega)$. Then $h\in \mathcal{H}^2_{\mathrm{loc}}(\Omega)$ and for each $\Omega_0\subset  \Omega$ there is a constant $c= c(\Omega_0)$ independent of $h$ and $w$ such that: 
\begin{equation}\label{Elliptic}\|h\|_{\mathcal{H}^2(\Omega_0)}\le c \left(\|h\|_{L^2(\Omega)}+\|w\|_{L^2(\Omega)}\right).
\end{equation}
\end{theorem}

This together with the fact that the magnetic potential $A_{\mathrm{magn}}$ is twice continuously differentiable on $G_\epsilon$ and having in mind that $\epsilon\le2\varepsilon$ estimates the right-hand side of (\ref{Gepsilon}) via $\|g\|_2$ as follows
\begin{gather*}
\int_{\mathbb B_\epsilon}|\nabla g|^2\,d x\,d y\\\le 2^{10/3} \varepsilon^{4/3} c^2\left(\frac{1}{4^{1/3}}+ 2^{1/3}\right)\mathrm{max}\left\{1, \frac{4 C(2)}{\sqrt{\pi}}\right\}\,
\int_\Omega\left(|(i\nabla+A_{\mathrm{magn}})^2g|^2+|g|^2\right)\,d x\,d y.
\end{gather*}

The above inequality together with the following lemma completes the proof: 

\begin{lemma}\label{Delta}[The proof will be given in the end of this section]
For any $z\in \mathrm{Dom}(H_\Omega(A_{\mathrm{magn}}))$ the following estimate is valid
$$\int_\Omega |(i\nabla+A_{\mathrm{magn}})^2z+ z|^2\,d x\,d y \ge \int_\Omega (|(i\nabla+A_{\mathrm{magn}})^2z|^2+ |z|^2)\,d x\,d y.$$
\end{lemma}

\qed

\subsection{Proof of Lemma \ref{magnetic1}}
\setcounter{equation}{0}

Using the same notations as in previous subsection,
Lemma \ref{marchenko} and repeating the similar calculations one infers that
\begin{eqnarray} \nonumber\int_{\mathbb B_\epsilon} |g|^2 \,d x\,d y\le 2\epsilon^{4/3} \left(\frac{1}{4^{1/3}}+ 2^{1/3}\right)\mathrm{max}\left\{1, \frac{4 C(2)}{\sqrt{\pi}}\right\}\,\|g\|_{\mathcal{H}^1(G_\epsilon)}^2\\\label{epsilon1} \le2\epsilon^{4/3} \left(\frac{1}{4^{1/3}}+ 2^{1/3}\right)\mathrm{max}\left\{1, \frac{4 C(2)}{\sqrt{\pi}}\right\}\,\|g\|_{\mathcal{H}^2(G_\epsilon)}^2.\end{eqnarray}

We finish the proof similarly as in previous subsection.

\qed

\subsection{Proof of Lemma \ref{auxiliary1}}

We first prove the following auxiliary statement: there exists $\tau\in (\varepsilon, 2\varepsilon)$ such that

\begin{equation}\label{statement}
\int_0^{2\pi} \left|\nabla g(\tau \cos \varphi, \tau \sin \varphi)\right|^2\,d \varphi \le \frac{1}{\varepsilon^2} 
\|g\|_1^2.
\end{equation}

Let us assume the opposite: for any $r\in (\varepsilon, 2\varepsilon)$ we have 
$$
\int_0^{2\pi} \left|\nabla g(r \cos \varphi, r \sin \varphi)\right|^2\,d \varphi>\frac{1}{\varepsilon^2} \|g\|_1^2.
$$

Passing to polar coordinates in integral $\int_{\mathbb B_{2\varepsilon}\backslash \mathbb B_\varepsilon} |\nabla g|^2\,d x\,d y$ and using the above bound we get
$$\int_{\mathbb B_{2\varepsilon}\backslash \mathbb B_\varepsilon} |\nabla g|^2\,d x\,d y= \int_\varepsilon^{2\varepsilon} \,\int_0^{2\pi} r  \left|\nabla g(r \cos \varphi, r \sin \varphi)\right|^2\,d \varphi\, d r>\|g\|_1^2.$$

This contradicts with the fact that the left-hand side of the above inequality does not exceed $\|g\|_1^2$. Hence there exists at least one number $\tau\in (\varepsilon, 2\varepsilon)$ such that 

\begin{equation}\label{varphi0}
\int_0^{2\pi} |\nabla g(\tau \cos \varphi, \tau \sin \varphi)|^2\,d \varphi\le \frac{1}{\varepsilon^2} \|g\|_1^2.
\end{equation}

In view of the representation of the derivative of function $\tilde{g}(r, \varphi)=g(r\cos\varphi, r\sin\varphi)$
$$\frac{\partial \widetilde{g}}{\partial\varphi}(\tau, \varphi)= -\tau \frac{\partial g}{\partial x}(\tau \cos \varphi, \tau \sin \varphi) \sin \varphi +\tau \frac{\partial g}{\partial y}(\tau \cos \varphi, \tau \sin \varphi) \cos \varphi$$ 
and the Cauchy inequality we conclude that 
$$\left|\frac{\partial \widetilde{g}}{\partial\varphi}(\tau, \varphi)\right|\le 2\varepsilon \left| \nabla g(\tau \cos \varphi, \tau \sin \varphi)\right|.
$$

Hence, employing (\ref{varphi0}) we establish that
$$
\int_0^{2\pi}\left|\frac{\partial \tilde{g}}{\partial\varphi}(\tau, \varphi)\right|^2\,d \varphi \le 4\varepsilon^2 \int_0^{2\pi}\left|\nabla g(\tau \cos \varphi, \tau \sin \varphi)\right|^2\,d \varphi\le 4\|g\|_1^2,
$$
which concludes the proof.\qed 

\subsection{Proof of Lemma\,\ref{omega.}}

Let $(x_0, y_0)\in \Omega_{\Gamma_\varepsilon}$. Assume the validity of the first condition of \textbf{property$^*$}. Let $y_1(x_0)$ be a point of intersection of $l(x_0)$ and the boundary of $\Omega$. Without less of generality suppose that $y_1(x_0)< y_0$ and the interval $(y_1(x_0), y_0)$ does not contain points from $\Gamma_\varepsilon$. 
One can easily check that there exists
$y_2(x_0)\in (y_1(x_0), y_0)$ such that
$$|u(x_0, y_0)|\le\frac{1}{\sqrt{y_0-y_1(x_0)}}\sqrt{\int_{y_1(x_0)}^{y_0}|u(x_0, z)|^2\,d z}.$$

Therefore 
\begin{eqnarray*}
|u(x_0, y_0)|^2= \left|u(x_0, y_2(x_0))+\int_{y_2(x_0)}^{y_0}\frac{\partial u}{\partial z}(x_0, z)\,d z\right|^2\\\nonumber\le 2|u(x_0, y_2(x_0))|^2+ 2(y_0-y_2(x_0))\int_{y_2(x)}^y\left|\frac{\partial u}{\partial z}(x, z)\right|^2\,d y\\\nonumber\le \frac{2}{y_0-y_1(x_0)}\int_{y_1(x_0)}^{y_0}|u(x_0, z)|^2\,d z
+ 2(y_0-y_2(x_0))\int_{y_2(x_0)}^{y_0}\left|\frac{\partial u}{\partial z}(x_0, z)\right|^2\,d z\\\le 
\frac{2}{\mathrm{dist}((x_0, y_0), \partial \Omega)}\int_{y_1(x_0)}^{y_0}|u(x_0, z)|^2\,d z
+ 2 \mathrm{diam}(\Omega)\int_{y_1(x_0)}^{y_0}\left|\frac{\partial u}{\partial z}(x_0, z)\right|^2\,d z\\\le\frac{2}{\mathrm{dist}((x_0, y_0), \partial \Omega)}\int_{l(x_0)\setminus\Gamma_\varepsilon}|u(x_0, z)|^2\,d z
+ 2 \mathrm{diam}(\Omega)\int_{l(x_0)\setminus\Gamma_\varepsilon}\left|\frac{\partial u}{\partial z}(x_0, z)\right|^2\,d z,
\end{eqnarray*} 
where  $\mathrm{diam}(\Omega)$ is the diameter of $\Omega$ and $\mathrm{dist}((x_0, y_0), \partial\Omega)$ is the distance between $(x_0, y_0)$ and the boundary of $\Omega$. This proves (\ref{curve1}) with $$C_8=\left(\mathrm{max}\left\{\frac{2}{\mathrm{dist}((x_0, y_0), \partial \Omega)},\, 2 \mathrm{diam}(\Omega)\right\}\right)^{/2}.$$

The case when $(x_0, y_0)$ satisfies second condition of \textbf{property$^*$} can be studied similarly. Repeating the same ideas one gets the validity of (\ref{curve2}).\qed

\subsection{Proof of Lemma\,\ref{Delta}}

It is straightforward to check that
\begin{eqnarray}\nonumber\int_\Omega |(i\nabla+A_{\mathrm{magn}})^2u+ u|^2\,d x\,d y\\\nonumber=
\int_\Omega |(i\nabla+A_{\mathrm{magn}})^2u|^2\,d x\,d y+2 \int_\Omega |i\nabla u+A_{\mathrm{magn}}u|^2\,d x\,d y+\int_\Omega |u|^2\,d x\,d y\\\label{resolvent}\ge\int_\Omega |u|^2\,d x\,d y.\end{eqnarray}

Let us consider two cases:
\begin{eqnarray}
\label{delta1}&&\int_\Omega |(i\nabla+A_{\mathrm{magn}})^2u|^2\,d x\,d y \ge 4\int_\Omega |u|^2\,d x\,d y,
\\\label{delta2}&&\int_\Omega |(i\nabla+A_{\mathrm{magn}})^2u|^2\,d x\,d y< 4\int_\Omega |u|^2\,d x\,d y.\end{eqnarray}
Starting from the first one and employing (\ref{resolvent}) we have
\begin{eqnarray*}&&\sqrt{\int_\Omega |(i\nabla+A_{\mathrm{magn}})^2u+u|^2\,d x\,d y}\ge \sqrt{\int_\Omega |(i\nabla+A_{\mathrm{magn}})^2u|^2\,d x\,d y}-
\sqrt{\int_\Omega |u|^2\,d x\,d y} \\&&\ge\frac{1}{2} \sqrt{\int_\Omega |(i\nabla+A_{\mathrm{magn}})^2u|^2\,d x\,d y}\ge  \frac{1}{4}\sqrt{\int_\Omega |(i\nabla+A_{\mathrm{magn}})^2u|^2\,d x\,d y}+\frac{1}{2}\sqrt{\int_\Omega |u|^2\,d x\,d y}\\&&\ge\frac{1}{4}\left(\sqrt{\int_\Omega |(i\nabla+A_{\mathrm{magn}})^2u|^2\,d x\,d y}+\sqrt{\int_\Omega |u|^2\,d x\,d y}\right).
\end{eqnarray*}
Hence we arrive at the bound
\begin{equation}\label{last}\int_\Omega |(i\nabla+A_{\mathrm{magn}})^2u+u|^2\,d x\,d y\ge \frac{1}{16}\left(\int_\Omega |(i\nabla+A_{\mathrm{magn}})^2u|^2\,d x\,d y+\int_\Omega |u|^2\,d x\,d y\right).
\end{equation}
Now let us consider the case (\ref{delta2}). In view of inequality (\ref{resolvent}) we conclude
\begin{eqnarray*}
\int_\Omega |(i\nabla+A_{\mathrm{magn}})^2u+u|^2\,d x\,d y\ge \int_\Omega |u|^2\,d x\,d y\ge\frac{1}{2}\int_\Omega |u|^2\,d x\,d y+ \frac{1}{8} \int_\Omega |(i\nabla+A_{\mathrm{magn}})^2u|^2\,d x\,d y\\\ge \frac{1}{8}\left(\int_\Omega |(i\nabla+A_{\mathrm{magn}})^2u|^2\,d x\,d y+ \int_\Omega |u|^2\,d x\,d y\right).
\end{eqnarray*}

Combining the above estimate together with (\ref{last}) we complete the proof of the lemma. \qed

%\subsection*{Acknowledgements}

%The work of D.B. is supported by the Czech Science Foundation (GA\v{C}R) within the project  21-07129S.

\bigskip

\section*{Funding} The research of D.B. and B.S. was supported by the Czech-Polish project BPI/PST/2021/1/00031. Furthermore, D.B. wants to thank for the support to the Czech Science Foundation (GACR) within the project 21-07129S.

\end{document}